\def\p{{\pi}} 
\def\s{{\sigma}}

\def\B{{\cal B}}

\def\F{{\cal F}}

\def\N{{\cal N}} 
\def\P{{\cal P}} 
 
\def\S{{\cal S}} 
\def\T{{\cal T}}

\def\ssf{{\mathsf f}} 
\def\sg{{\mathsf g}} 
\def\sk{{\mathsf k}}

\def\pf{{\hfill$\Box$}} 
\def\proof{{\noindent {\it Proof}.\ \ }} 
\def\ra{{\rightarrow}}

\def\l({{\left(}} 
\def\r){{\right)}}

\def\({{\Biggl(}} 
\def\){{\Biggr)}} 
\def\[{{\Biggl[}} 
\def\]{{\Biggr]}}

\documentclass[12pt]{article}

\newtheorem{thm}{Theorem}

\usepackage{graphicx,amssymb,amsmath,verbatim} 

\begin{document}

%
%
\title{Two New Bijections on Lattice Paths} 

\author{ 
Glenn Hurlbert\thanks{ 
    \texttt{hurlbert@asu.edu}, corresponding author 
    }\\ 
Vikram Kamat\thanks{ 
    \texttt{vikram.kamat@asu.edu} 
    }\\ 
 \\ 
Department of Mathematics and Statistics\\ 
Arizona State University\\ 
Tempe, Arizona 85287-1804 
} 
\maketitle 


\vspace{0.5 in} 

%
%
\begin{abstract} 
Suppose $2n$ voters vote sequentially for one of two candidates.  
For how many such sequences does one candidate have strictly more votes 
than the other at each stage of the voting? 
The answer is $\binom{2n}{n}$ and, while 
easy enough to prove using generating functions, for example, 
only two combinatorial proofs exist, due to Kleitman and Gessel.
In this paper we present two new (far simpler) bijective proofs.  
\vspace{0.1 in}

\end{abstract} 

\newpage

%
%
\section{Introduction}\label{Intro} 

Suppose $A$ and $B$ are candidates and there are $2n$ 
voters, voting sequentially.  
In how many ways can $A$ get $n$ votes and $B$ get $n$ votes 
such that $A$ is always ahead of or tied with $B$? This is the 
famous ballot problem, the solution of which is counted by the 
Catalan number $C_n.$ A more general version of this problem is stated 
in terms of probability: let candidate $A$ receive $a$ votes and $B$ 
receive $b$ votes and compute the probability that $A$ never falls 
behind $B$. This is known to be $\frac{a-b}{a+b}$, first proved by
Andr\'e \cite{A} using his reflection principle (this also appears in
\cite{F} and again in \cite{HP}).
A $q$-binomial variation of the problem appears in \cite{KM}, and a weighted
variation is found in \cite{EFY}.
Some have also considered an $n$-dimensional version by 
generalizing Andr\'e's reflection proof to the many candidate ballot 
problem -- one of these approaches can be found in \cite{Z}. 
Others \cite{FR, Lo,Ly} have considered the situation in which the
number of votes of the two candidates remains close.
In this paper, we concern ourselves with another variation 
of the two candidate ballot problem discounting all instances when 
the two candidates are tied. It can be more formally expressed in 
terms of plus minus sequences.  

Let $\S_n=\{s_1\ldots s_{2n}\ |\ s_i\in\{-1,+1\}\}$.  
For $S\in\S_n$ let $\s(S)=(\s_1,\ldots,\s_{2n})$, where $\s_i=\sum_{j=1}^is_i$.  
We write that $\s\not=0$ (resp, $\s>0$, $\s<0$) when each $\s_i\not=0$ (resp.  
$\s_i>0$, $\s_i<0$), and call $S$ {\it zero-free} (resp. {\it positive}, 
{\it negative}) if $\s(S)\not=0$ (resp. $\s(S)>0$, $\s(S)<0$).  
The set of zero-free (resp. positive, negative) sequences of $\S_n$ is 
denoted by $\F_n$ (resp. $\P_n$, $\N_n$).  
We will find it useful to denote $\s_{2n}$ by $\sum S$.
A sequence $S\in\S_n$ is {\it balanced} if $\sum S=0$, and we use $\B_n$ 
to denote all balanced sequences in $\S_n$.  
We also use the notations $\S_n^+, \F_n^+, \B_n^+$ to denote those sequences 
that start with $+1$ (with the obvious analogous definitions for 
$\S_n^-$, etc.).  
Note that $\F_n^+=\P_n$, $\F_n^-=\N_n$, and $|\P_n|=|\N_n|$.  

It is known that $|\F_n|=|\B_n|$ for every $n$ and, to our knowledge,
only two bijections have appeared in print, due to Kleitman \cite{K} and
Gessel (see \cite{S}).\footnote{It has come to our attention recently that
a combinatorial proof similar to ours is given by Callan \cite{C}.}
Here we give two new bijections for this result, one indirect 
(Section \ref{KP}) and one direct (Section \ref{Hurl}).  
Of course, it is enough to show these for sequences in $\S_n^+$.  
That our direct bijection $\sg$ (see Section \ref{Hurl}) 
differs from that of Kleitman's bijection $\sk$ can be seen for 
virtually any sequence $P \in P_n$ for any $n$, such as the 
following.  
\begin{eqnarray*} 
  P & = & ++-+++---+++--++++-+--+++++- \\ 
  \sg(P) & = & +-+---++++--++++--+-++-----+ \\ 
  \sk(P) & = & ---+++---+++----++-+--+++++- \\ 
\end{eqnarray*} 
In addition, ours is considerably simpler to navigate.  
Like Kleitman's, Gessel's bijection does not preserve the first coordinate.

As is well known, one of the interesting applications of this result is the 
derivation of the generating function $F(x)$ for the sequence 
$\{\binom{2n}{n}\}$.  
Indeed, the factoring of a sequence $S\in\S_n$ into its maximum length 
balanced initial subsequence and corresponding terminal zero-free subsequence 
results in 
the relation $|\S_n|=\sum_k|\B_k||\F_{n-k}|$, from which the bijection 
yields 
\begin{equation}\label{sved}
4^n=\sum_k\binom{2k}{k}\binom{2n-2k}{n-k}\ ,
\end{equation}
one short convolutional 
step from proving that $F(x)=(1-4x)^{-1/2}$.  
Marta Sved \cite{S} recounts the history of Identity (\ref{sved})
and notes some combinatorial proofs of it that were submitted by readers. 

Most satissfyingly, this paper is the result of the first author challenging 
the students in his graduate combinatorics course to find such a bijection.  

In this article we use the notations $[s]=\{1,2,\ldots,s\}$ and 
$(r,s)=[r+1,s-1]=\{r+1,\ldots,s-1\}$.

%
%
\section{Proofs}\label{Proofs} 

%
%
\subsection{Indirect Bijection}\label{KP}

Let $\S_{n,k}^+$ be the set of all $S\in\S_n^+$ with 
$\sum S=2k$. Such an $S$ contains $(n-k)$ `$-1$'s and $(n+k)$ 
`$1$'s. Therefore, $$|\S_{n,k}^+|=\binom{2n-1}{n-k}\ ,$$ 
since each such $S$ begins with a `1'.  
For each $1\leq k\leq n,$ let $\P_{n,k}$ be the set of all $P\in 
\P_n$ having $\sum P=2k$.  
Then $|\P_n|=\sum_{k=1}^n|\P_{n,k}|$.  
Define $\T_{n,k}^+=S_{n,k}^+-\P_{n,k}$, so that 
$$|\P_{n,k}|=|\S_{n,k}^+|-|\T_{n,k}^+|=\binom{2n-1}{n-k}-|\T_{n,k}^+|\ .$$ 
Finally, let $\S_{n,k}^-$ be all sets $S\in S_n^-$ for which $\sum S=2k$.  
Because such an $S$ begins with a `$-1$', we have 
$$|\S_{n,k}^-|=\binom{2n-1}{n-(k+1)}\ .$$ 
If $|\T_{n,k}^+|=|\S_{n,k}^-|$, then 
\begin{eqnarray*} 
\sum_{k=1}^{n}|\P_{n,k}| 
&=&\sum_{k=1}^{n}\left[\binom{2n-1}{n-k}-\binom{2n-1}{n-(k+1)}\right]\\ 
&=&\binom{2n-1}{n-1}\ , 
\end{eqnarray*} 
by telescoping.  
Thus, we obtain $|\F_n|=\binom{2n}{n}=|\B_n|$ if 
$|\T_{n,k}^+|=|\S_{n,k}^-|$. We show that this indeed holds by 
demonstrating a bijection between the two sets.  
\medskip 

\noindent 
{\bf From $\T_{n,k}^+$ to $\S_{n,k}^-$.} 
\smallskip 

Let $T\in \T_{n,k}^+$. We can factor $T$ as follows: 
$$T=T_1T_2$$ where $T_1$ is the smallest balanced initial subsequence, 
say of length $m$.  
(In lattice path language, $m$ is the first step that hits the diagonal 
$y=x$.) 
By definition every sequence in $\T_{n,k}^+$ can be 
represented this way. Let $\overline{T_1}$ be the subsequence of 
length $m$ obtained by negating every term of $T_1.$ Now, define a 
function $f \colon \T_{n,k}^+ \to \S_{n,k}^-$ as follows: 
$$f(T)=f(T_1T_2)=\overline{T_1}T_2\ .$$ 

Consider an example. Let 
\begin{eqnarray*} 
T&=&++-++-+--\fbox{$-$}++-++++++-+\ .  
\end{eqnarray*} 
Here, we see that $T\in \T_{11,4}^+$ and $m=10$ (with the $10^{\rm th}$ 
position boxed).  
Also, 
\begin{eqnarray*} 
f(T)&=&--+--+-++\fbox{+}++-++++++-+\ .  
\end{eqnarray*} 

Clearly, $f(T)$ belongs to $\S_{11,4}^-.$ The reason why this will 
hold in general is because only a balanced subsequence is negated 
and so $\s$ is unchanged.  
Moreover, the 
smallest balanced subsequence always contains the first 
element(which is a $1$) and thus, $f(T)$ will always have $-1$ as 
its first element.  
\medskip 

\noindent 
{\bf From $\S_{n,k}^-$ to $\T_{n,k}^+$.} 
\smallskip 

Analogously, define a function $g \colon \S_{n,k}^- \to \T_{n,k}^+$ 
as follows: for any $S\in \S_{n,k}^-$, 
$$g(S)=g(S_1S_2)=\overline{S_1}S_2\ ,$$ 
where $S=S_1S_2$ with $S_1$ its smallest balanced initial subsequence, and 
$\overline{S_1}$ the negation of $S_1$.  

In the example above, we have 
\begin{eqnarray*} 
S&=&--+--+-++\fbox{+}++-++++++-+\ \in\  \S_{11,4} 
\end{eqnarray*} 
and 
\begin{eqnarray*} 
g(S)&=&++-++-+--\fbox{$-$}++-++++++-+\ .  
\end{eqnarray*} 

Clearly, $g(S)\in \T_{11,4}^+.$  As above, this holds in general since 
$\s$ is unchanged and $g(S)$ begins with a `1' because $S$ begins with a 
`$-1$'.  

\begin{thm} 
The functions $f \colon \T_{n,k}^+ \to \S_{n,k}^-$ and $g \colon 
\S_{n,k}^- \to \T_{n,k}^+$ are inverses of each other.  
\end{thm} 

\proof 
Let $T\in \T_{n,k}^+$ and let $S=f(T).$ We will show that $g(S)=T.$ 
We have $T=T_1T_2$ where $|T_1|=m$, as defined above, and 
$S=\overline{T_1}T_2.$ Now, when we apply function $g$ to $S$, since 
$T_1$ is the smallest balanced initial subsequence in $T$, $\overline{T_1}$ 
is the smallest balanced initial subsequence in $S$. This gives $g(S)=T.$ 

A similar argument proves that, for any $S\in\S_{n,k}^-$, 
$f(T)=S$ when $g(S)=T$.  
\pf 


%
%
\subsection{Direct Bijection}\label{Hurl} 

As previously mentioned, we will give a bijection between $\B_n^+$ and $\P_n$.  
The obvious bijection between $\B_n^-$ and $\N_n$, and hence $\B_n$ and $\F_n$ 
follows.  
\medskip 

\noindent 
{\bf From $\B_n^+$ to $\P_n$.} 
\smallskip 
For $B\in\B_n^+$ we define a set $\p(B)$ of {\it peaks} of $B$ as follows.  
Set $t=\max\{\s_i(B)\}$ and for $1\le k\le t$ 
let $\p_k$ be the index of the left-most occurrence of $k$ in $\s=\s(B)$: 
$\p_k=\min\{i\ |\ \s_i=k\}$.  
For example, if 
$$B\ =\ \fbox{+}-+---+++\fbox{+}--++\fbox{+}\fbox{+}--+-++-----+$$ 
then 
{\footnotesize
$$\s(B)\ =\ (\fbox{1},0,1,0,-1,-2,-1,0,1,\fbox{2},1,0,1,2,\fbox{3},\fbox{4},3,2,3,2,3,4,3,2,1,0,-1,0)$$}
and 
$$\p(B)=\{1,10,15,16\}\ .$$ 
We have boxed in  the peak locations as shown.  

Next we define the set of intervals $I_k=(\p_k,\p_{k+1})$, with 
$I=\cup I_k=[2n]-\p(B)$.  
(Artificially, we set $\p_{t+1}=2n+1$ in order to define 
$I_t$; in this case we have $\p_5=29$ and $I_4=(16,29)=[17,28]$.) 
The key property here is that $\s_i\le k$ for every $i\in I_k$.  

Finally we define $\ssf=\ssf(B)$ by $\ssf_i=B_i$ for all $i\in\p(B)$ and $\ssf_i=-B_i$ 
otherwise (for all $i\in I$).  
For this example, we obtain 
$$\ssf\ =\ \fbox{+}+-+++---\fbox{+}++--\fbox{+}\fbox{+}++-+--+++++-$$ 
and 
{\small
$$\s(\ssf)\ =\ (\fbox{1},2,1,2,3,4,3,2,1,\fbox{2},3,4,3,2,\fbox{3},\fbox{4},5,6,5,6,5,4,5,6,7,8,9,8)\ .$$} 
Note that we have $\s(\ssf)>0$, so that $\ssf(B)\in\P_n$.  
This holds in general for the following reason.  
By the definition of $\p_k$ we have $\sum B_{\p_{k+1}-1}=\sum B_{\p_k}$.  
This means that $B$ is balanced on each interval $I_k$ with $k<t$.  
Hence each $\sum\ssf_{\p_k}=\sum B_{\p_k}$ and thus $\s(\ssf)_i\ge k$ 
for every $i\in I_k$.  
In particular, $\ssf_i\ge 1$ for all $i$.  
\medskip 

\noindent 
{\bf From $\P_n$ to $\B_n^+$.} 
\smallskip 

For $P\in\P_n$ we define a set $\Pi(P)$ of {\it pivots} of $P$ as follows.  
Set $T=\frac{1}{2}\sum P$, let $\Pi_1=1$, and for $1<k\le T$ let 
$\Pi_k$ be one more than the index of the right-most occurrence of $k-1$ in 
$\s=\s(P)$: $\Pi_k=1+\max\{j\ |\ \s_j=k-1\}$.  
For example, if 
$$P\ =\ \fbox{+}+-+++---\fbox{+}++--\fbox{+}\fbox{+}++-+--+++++-$$ 
then 
{\small
$$\s(P)\ =\ (\fbox{1},2,1,2,3,4,3,2,1,\fbox{2},3,4,3,2,\fbox{3},\fbox{4},5,6,5,6,5,4,5,6,7,8,9,8)$$} 
and 
$$\Pi(P)=\{1,10,15,16\}\ .$$ 
We have boxed in  the pivot locations as shown.  

Next we define the set of intervals $J_k=(\Pi_k,\Pi_{k+1})$, with 
$J=\cup J_k=[2n]-\Pi(P)$.  
(Artificially, we set $\Pi_{T+1}=2n+1$ in order to define 
$J_T$; in this case we have $\Pi_5=29$ and $J_4=(16,29)=[17,28]$.) 
The key property here is that $\s_j\ge k$ for every $j\in J_k$.  
Finally we define $\sg=\sg(P)$ by $\sg_j=P_j$ for all $j\in\Pi(P)$ and $\sg_j=-P_j$
otherwise (for all $j\in J$).
For this example, we obtain
$$\sg\ =\ \fbox{+}-+---+++\fbox{+}--++\fbox{+}\fbox{+}--+-++-----+$$
and
{\footnotesize
$$\s(\sg)\ =\ (\fbox{1},0,1,0,-1,-2,-1,0,1,\fbox{2},1,0,1,2,\fbox{3},\fbox{4},3,2,3,2,3,4,3,2,1,0,-1,0)\ .$$}
Note that we have $\sum\sg=0$, so that $\sg(P)\in\B_n^+$.
This holds in general for the following reason.
By the definition of $\Pi_k$ we have $\sum P_{\Pi_{k+1}-1}=\sum P_{\Pi_k}$.
This means that $P$ is balanced on each interval $J_k$ with $k<T$.
Hence each $\sum\sg_{\Pi_k}=\sum P_{\Pi_k}$ and thus $\s(\sg)_j\le k$
for every $j\in J_k$.
In particular, by the definition of $T$ we have $\sum_{j\in J_T}P_j=T$,
so that $\sum_{j\in J_T}\sg_j=-T$, and hence $\sum \sg=0$.

\begin{thm}
The functions $\ssf:\B_n^+\ra\P_n$ and $\sg:\P_n\ra\B_n^+$ are
bijections and, in fact, inverses of each other.
\end{thm}

\proof
The arguments above show that $\ssf$ and $\sg$ are well-defined.
That they are bijections will follow from their inverse relationship.

We suppose first that $\ssf(B)=P$ for $B\in\B_n^+$, and show that $\sg(P)=B$.
This will follow from showing inductively that $\Pi(P)=\p(B)$.
Of course, $\Pi_1=1=\p_1$, so assume that $\Pi_k=\p_k$.
Because we know for all $k$ that $\s(B)_i\le k$ for all $i\in I_k$
and that $\s(B)_{\p_{k+1}}=k+1$, we therefore know for all $k$ that
$\s(P)_i\ge k$ for all $i\in I_k$ and that $\s(P)_{\p_{k+1}}=k+1$.
In particular, the right-most occurrence of $k$ in $\s(P)$ occurs with
index $\p_{k+1}-1$; i.e. $\Pi_{k+1}=\p_{k+1}$.

Next we suppose that $\sg(P)=B$ for $P\in\P_n$, and show that $\ssf(B)=P$.
As above, we show that $\p(B)=\Pi(P)$ by induction.
Again, $\p_1=1=\Pi_1$, so we assume that $\p_k=\Pi_k$.
We know for all $k$ that $\s(P)_j\ge k$ for all $j\in J_k$
and that $\s(P)_{\Pi_{k+1}}=k+1$, and thus we know for all $k$ that
$\s(B)_j\le k$ for all $j\in J_k$ and that $\s(B)_{\Pi_{k+1}}=k+1$.
In particular, the left-most occurrence of $k+1$ in $\s(B)$ occurs with
index $\Pi_{k+1}$; i.e. $\p_{k+1}=\Pi_{k+1}$.
\pf

%
%


%
%

\section*{Acknowledgement}
We thank Ira Gessel and David Callan for several helpful comments.

%
%
\bibliographystyle{plain}
%

%
%

\end{document}